\documentclass[12pt]{amsart}
\usepackage{graphicx}

\usepackage[pagewise]{lineno}

\newcommand{\Frac}[2]{{\textstyle\frac{#1}{#2}}}
\newcommand{\pfrac}[2]{\left(\frac{#1}{#2}\right)}

\begin{document}
\setcounter{page}{1}

\title{Evaluation of the Gauss integral}
\author[Dimitri Martila, Stefan Groote]{Dmitri Martila, Stefan Groote\\[12pt]
University of Tartu, Institute of Physics, W.~Ostwaldi~1, 50411~Tartu, Estonia}
\thanks{eestidima@gmail.com}

\begin{abstract}
The normal or Gaussian distribution plays a prominent role in almost all
fields of science. However, it is well known that the Gauss (or
Euler--Poisson) integral over a finite boundary, as it is necessary for
instance for the error function or the cumulative distribution of the normal
distribution, cannot be expressed by analytic functions. This is proven by the
Risch algorithm. Still, there are proposals for approximate solutions. In this
paper, we give a new solution in terms of normal distributions by applying a
geometric procedure iteratively to the problem.\\[12pt]
MSC Class: 62E17; 60E15; 26D15.
\end{abstract}

\maketitle

\section{Introduction}
The normal or Gaussian distribution plays a prominent role in almost all
fields of science, as the sum of random variables tends to the normal
distribution if the quite general conditions of the central limit
theorem~\cite{Kendall:1965} are satisfied. Besides the unbounded normal
integral, the bounded integral or error function is crucial for the
determination of probabilities. However, there is no analytic expression
found for this function, a fact that can be tested by using the Risch
algorithm~\cite{Risch:1969,Risch:1970}. Powerful modern computer facilities
but also simple personal computers allow for a numerical calculation of the
error function with any needed precision. However, if a multitude of such
calculations has to be performed in a limited time, for instance, in Monte
Carlo simulations, the processing time becomes essential. In order to speed up
these calculations, simple and more educated approximations have been proposed
in the literature. The spectrum of approximations contains, for instance, the
Gaussian exponential function, including either numerical
constants~\cite{Ordaz:1991} or powers and square roots~\cite{Shenton:1954},
approximations using ordinary and hyperbolic tangent
functions~\cite{Vazques-Leal:2012}, a rational function of an exponential
function with the exponent given by a power
series~\cite{Sandoval-Hernandez:2019}, or an approximation by Jacobi theta
functions~\cite{Zhang:2018}. Without knowing the error function explicitly,
expectation values can be calculated by an approximation of the normal
distribution by a series in ordinary exponential
functions~\cite{Chesnaeu:2019}.

The present paper contains a continuation of this topic. Employing a geometric
approach, we provide an approximation of the squared error function by a
finite sum of $N$ Gaussian exponential functions with different widths, where
the values of which are constrained to fixed intervals. We show that, by
fine-tuning these width parameters, one can optimise the precision, which,
even for the leading order $N=1$, is better than the error estimates given by
the constraints in Ref.~\cite{Shenton:1954}. In addition, by choosing $N$ as
appropriately large, one can afford an arbitrary precision. On the other hand,
even on a personal computer, the calculation with our leading order
approximation is obtained 34 times faster than an exact numerical calculation,
the processing time for higher orders being multiplied by $N$.

Our paper is organised as follows. In Sec.~2 we introduce the basic concepts
for the calculation of the Gaussian integral that are necessary for the
understanding of our geometric approach. The precisions of the leading order
approximation obtained here and simple, straightforward extensions of this
approximation are discussed in Sec.~3. In Sec.~4 we explain the geometric
background for our approximation and provide a systematic way to create
higher order approximations. The iterative construction of higher order
approximations is explained in general in Sec.~5 in terms of partitions before
we turn to the partition into $N=2^p$ intervals for increasing values of $p$.
In Sec.~6 we explain a similar ternary construction. In Sec.~7 we provide our
conclusions and an outlook on possible extensions. The convergence of our
iterative procedure is discussed in more detail in Appendix~A. In addition, we
discuss the continuum limit, which, of course, cannot be part of the algorithm
but allows, as a bonus, for a different representation of the error function. 

\section{Basic concepts}
\begin{figure}[ht]
\includegraphics[width=0.95\linewidth]{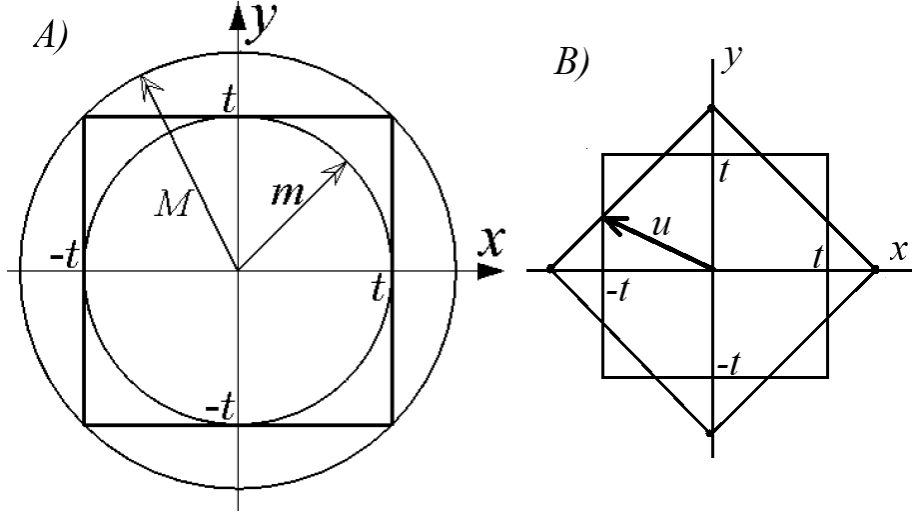}
\caption{(A) As the integrant is positive, the value of the integral over the
  interior of the square lies between the values of integrals over the
  interiors of the circumferences of radii $m$ and $M$. (B) the previous
  integration square is taken and rotated by $\pi/4$ and put together with the
  exact copy of the square. The integral over the common inner area is denoted
  by $\omega$.\label{fig1}}
\end{figure}

The error function is based on the standard normal density distribution
\begin{equation}
\rho(x)=\frac{1}{\sqrt{2\pi}}e^{-x^2/2}
\label{Eq1}
\end{equation}
which does not have a direct practical meaning, while it is desirable to
evaluate the integral of this function over a bounded interval $[-t,t]$,
leading to the probability $P(t)$ to find the result within this interval,
\begin{equation}
P(t)=\int_{-t}^t\rho(x)\,dx=\int_{-t}^t\rho(y)\,dy\,.
\label{Eq2}
\end{equation}
From Eqs.~(\ref{Eq1}) and~(\ref{Eq2}), the square of probability is given by
\begin{equation}
P^2(t)=\frac{1}{2\pi}\,\int_{-t}^t\rho(x)\,dx\,\int_{-t}^t\rho(y)\,dy=
\frac{1}{2\pi}\,\int_{-t}^t\int_{-t}^te^{-(x^2+y^2)/2}\,dx\,dy\,,
\end{equation}
where the integration area is a square in Fig.~\ref{fig1}(A). Introducing
polar coordinates $x = r\,\cos\,\varphi$ and $y = r\,\sin\,\varphi$, one
obtains
\begin{equation}
P^2=\frac{1}{2\pi}\,\int\int e^{-r^2/2}\,r\,dr\,d\varphi\,.
\label{Eq3}
\end{equation}
The integral in Eq.~(\ref{Eq3}) is analytically calculable if the integration
is performed over the interior of some circle with radius $R$. Indeed,
\begin{equation}
I^2(R)=\frac{1}{2\pi}\,\int^{2\pi}_0 d\varphi
  \int^R_0 e^{-r^2/2}\,r\,dr=1-e^{-R^2/2}\,.
\end{equation}
Here the function $I(R)$ increases monotonically with $R$ as the integral in
Eq.~(\ref{Eq3}) is taken over a positive function. This is why
$I(R = m) < P < I(R = M)$, with $m = t$ and $M = t\,\sqrt{2}$, see
Fig.~\ref{fig1}(A). Therefore,
\begin{equation}
P(t)=\sqrt{1-e^{-k^2(t)t^2/2}}\,,
\label{Eq4}
\end{equation}
where $1<k(t)<\sqrt 2$. Using a PC for analyzing the set of
Eqs.~(\ref{Eq1}), (\ref{Eq2}), and~(\ref{Eq4}), one concludes that $k(t)$
is even more constrained by $1<k(t)<\sqrt{4/\pi}$. Hence,
\begin{equation}
P_m(t)<P(t)<P_M(t),
\label{Eq44}
\end{equation}
and $P_M(t)-P_m(t)$ has a maximum of $0.0592$ at $t=t_0=1.0668$. The
inequality~(\ref{Eq44}) proves to be incomparably elegant, easy to remember,
and much more accurate than the best result of Ref.~\cite{Shenton:1954},
which, if transformed into the present formalism, will be
\begin{equation}
P_m(t)=1-4\sqrt{\frac{2}{\pi}}\,\frac{{\rm exp}(-t^2/2)}{3\,t+\sqrt{t^2+8}}\,,
\end{equation}
\begin{equation}
P_M(t)=1-\frac{1}{\sqrt{2\pi}}\,(\sqrt{t^2+4}-t)\,{\rm exp}(-t^2/2)\,.
\end{equation}
The largest range $P_M(t)-P_m(t)\approx 0.330$ for these constraints occurs
for $t = 0$. Compared to this, even at the leading order observed so far, our
value for $P_M(t)-P_m(t)<0.0592\ll 0.330$ is more restrictive. In more detail,
if, in Eq.~(\ref{Eq4}), we choose $k = \sqrt{4/\pi}$, the error will be below
$0.006$, but if we take $k = 1.116$, the maximal error is only $0.0033$. Even
modern reviews on this subject do not have better results~\cite{Latala:2002}.

\section{Simple extensions}
By adding additional terms to the leading order approximation, one can
increase the precision further. For the normal integral
\begin{equation}
P(t)=\frac{1}{\sqrt{2\pi}}\,\int^t_{-t}e^{-x^2/2}dx
\label{pppte}
\end{equation}
and the leading order approximation, for $k_1=1.116$, one has
\begin{equation}
\Big|P(t)-\sqrt{1-e^{-k_1^2t^2/2}}\Big|<0.0033\,.
\label{rteq}
\end{equation}
However, the precision increases by a factor of $14\approx 0.0033/0.00024$
by using
\begin{equation}
\Big|P(t)-\sqrt{1-\Frac12(e^{-k_1^2t^2/2}+e^{-k_2^2t^2/2})}\Big|<0.00024\,,
\label{rteq2}
\end{equation}
where $k_1=1.01$, $k_2=1.23345$. The next order of precision has
\begin{equation}\label{rteq3}
\Big|P(t)-\sqrt{1-\Frac13(e^{-k_1^2t^2/2}+e^{-k_2^2t^2/2}+e^{-k_3^2t^2/2})}
  \Big|<0.00003\,,
\end{equation}
where $k_1=1.02335$, $k_2=1.05674$, and $k_3=1.28633$. Therefore, this formula
with three exponentials is at least $8$ times more precise than the one with
two exponentials, and it is at least 110 times more precise than
Eq.~(\ref{rteq}) with one exponential only. Finally, it is
$11\,000\approx 0.33/0.00003$ times more precise than the approximation in
Ref.~\cite{Shenton:1954}. As it turns out, the values for the parameters
$k_i$ for $i$ running from $1$ to $N$ take values between $1$ and $\sqrt 2$,
while the sum of the exponential factors is divided by $N$. Still, there is a
degree of arbitrariness in the determination of these parameters. In order to
remove this arbitrariness, in the following, we develop an iterative method
based on geometry.

\section{Geometric background of our procedure}
In order to understand our method, we refer to Fig.~\ref{fig1}(B) for the
first step. Starting with the square with side length $2t$ representing the
square $P^2(t)$ of the probability, we turn this square by an angle of $\pi/4$
to obtain $P^2(t)$ again. In union and intersection, these two overlayed
squares construct two 8-angle figures. In order to obtain the area of the
larger figure, one has to subtract the area of the smaller figure from the
twofold square area, as this smaller figure is covered twice by the two
squares. Accordingly, for the integrals over the probability density, one
obtains the relation
\begin{equation}
\Omega(t)=2P^2(t)-\omega(t)
\label{2P2}
\end{equation}
between the probabilities. Now,
\begin{equation}
\omega(t)=1-e^{-k_1^2t^2/2}\,,\quad \Omega(t)=1-e^{-k_2^2t^2/2}\,,
\label{lay25}
\end{equation}
where
\begin{equation}
1<k_1<1/\cos\theta\,,\quad 1/\cos\theta<k_2<\sqrt{2},
\label{lay26}
\end{equation}
and the angle $\theta=\pi/8$ is enclosed between the $x$ axis and the vector
$u$ shown in Fig.~\ref{fig1}(B). We study $\Delta(k_1,k_2,t)=P(k_1,k_2,t)-P(t)$
with
\begin{equation}
P(k_1,k_2,t)=\sqrt{1-\Frac12\left(e^{-k_1^2t^2/2}+e^{-k_2^2t^2/2}\right)}\,.
\end{equation}
Drawing three-dimensional graphics and looking for a minimum of
$|\Delta(k_1,k_2,t)|$, one obtains
\begin{equation}
\Big|\Delta(k_1,k_2,t)\Big|<0.00024\,,
\end{equation}
where $k_1=1.01$ and $k_2=1.23345$. This is the starting point.

\section{Basic construction of the procedure}
In order to construct the iteration, one performs a partition of the figure
describing $\omega(t)$, $\Omega(t)$, or both of these, by repeating the
geometric construction shown before. For instance, taking only the larger
8-angle figure describing $\Omega(t)$, one can turn this figure by an angle
$\theta=\pi/16$ and overlay the new figure with the old one. In doing so, one
can separate a new larger and smaller 16-angle figure in the same way as
was carried out before for the 8-angle figures. Accordingly, by geometric
means, one obtains new constraints. In order to describe the procedure in a
unique way, in each iterative step, we rename $k_n$ by $k_{2n}$, and, if this
new $k_{2n}$ is subject to a partition, the smaller and larger figure of this
partition are related to the values $k_{2n-1}$ and $k_{2n}$, respectively.

Using the case in the previous section as an illustrative example for the
procedure, we might keep the smaller 8-angle figure related to $\omega(t)$ but
apply a partition to the larger 8-angle figure related to $\Omega(t)$.
Accordingly, $k_1$ is replaced by $k_2$ and $k_2$ is replaced by $k_4$, but this
new $k_4$ is again split up into $k_3$ and $k_4$. The constraint for the
lowest parameter $k_2$ (the former $k_1$) remains the same,
\begin{equation}
1\le k_2\le 1/\cos(\pi/8)
\end{equation}
whereas, for the two new higher parameters, we obtain
\begin{equation}
1/\cos(2\pi/16)\le k_3\le 1/\cos(3\pi/16)
\end{equation}
and
\begin{equation}
1/\cos(3\pi/16)\le k_4\le 1/\cos(4\pi/16)=\sqrt2.
\end{equation}
The intervals are consecutive, but $\pi/8$ is replaced by $2\pi/16$ in order
to indicate the new partition. Finally, the upper limit stays at
$1/\cos(4\pi/16)=1/\cos(\pi/4)=\sqrt 2$. For these values $k_2$, $k_3$, and
$k_4$, one obtains the approximation
\begin{equation}
P(t)\approx P(k_2,k_3,k_4,t)
\end{equation}
with
\begin{equation}
P(k_2,k_3,k_4,t)=\sqrt{1-\Frac12e^{-k_2^2t^2/2}-\Frac14e^{-k_3^2t^2/2}
  -\Frac14e^{-k_4^2t^2/2}}
\label{k1k2}
\end{equation}
because of the geometric transformations of Fig.~\ref{fig1} and the
corresponding double use of Eq.~(\ref{2P2}). Note that the set of parameters
$k_2$, $k_3$, $k_4$ is different from the set $k_1$, $k_2$ and $k_3$ in
Eq.~(\ref{rteq3}). Indeed, if, for Eq.~(\ref{k1k2}), one uses $k_2=1.025187$,
$k_3=1.1249$, and $k_4=1.31336$, the precision improves to $0.000015$.
Still, it is obvious that this example is only half of an iteration step, and
one could do definitely achieve a greater result by also performing the
partition for the smaller 8-angle figure, leading to four parameters separated
uniformly,
\begin{equation}
1\le k_1\le\frac1{\cos(\pi/16)}\le k_2\le\frac1{\cos(\pi/8)}\le k_3\le
\frac1{\cos(3\pi/16)}\le k_4\le\sqrt 2.
\label{k1k2k3k4}
\end{equation}
Therefore, a full iteration step is increasing the number of parameters $k_n$
by a factor of two, and, after $p$ full iteration steps, one has $N=2^p$
parameters. Each iteration step is finalized by optimizing the $N$ (or less)
parameters $k_n$. For any finite (or even very large) $N$, the constraints
\begin{equation}
k_n^{\rm min}(N)\le k_n\le k_n^{\rm max}(N)
\label{qnQ}
\end{equation}
with $k_n^{\rm min}(N)=k_{n-1}^{\rm max}(N)$ can be calculated from geometry
observations in a similar fashion. In practice, for a small set of parameters,
we use a graphical method. For instance, the method applied to obtain the
three values $k_2$, $k_3$, and $k_4$ in Eq.~(\ref{k1k2}) was to look for the
solution of the system of three equations
\begin{equation}\label{Qt0}
Q(t=1)=0\,,\quad Q(t=\sqrt{2})=0\,,\quad Q(t=2)=0\,,
\end{equation}
where
\begin{equation}
Q(t)=P^2(k_2,k_3,k_4,t)-P^2(t)\,.
\end{equation}
The values $t=1$, $\sqrt{2}$, and $2$ are used as nodes for this approximation.
Their choice depends on the application of the approximation and has to be
adjusted to the number of width parameters to be determined. Each equation
in~(\ref{Qt0}) can be treated individually. Therefore, the solution is very
easy to find. From $Q(t=1)=0$, one extracts the function $k_2=k_2(k_3,k_4)$.
Inserting these solutions into $Q(t=\sqrt{2})=0$, one extracts the two positive
functions $k_3=k_3(k_4)$ and $k_3=\tilde k_3(k_4)$. Inserting these solutions
into $Q(t=2)=0$ and plotting the function $Q(t=2,k_4)$, one finds the position
of the zero, which proves to be $k_4=1.31336$. Using this knowledge, one obtains
$k_3=k_3(k_4)$ and $k_2=k_2(k_3,k_4)$ as well. However, as $k_3=\tilde k_3(k_4)$
is given for $k_4<1$ only, this is not a valid solution, as $k_n\ge 1$ for all
$n$. Note that the graphical method cannot be applied any more for $N\ge 4$.
Instead, we used a random number generator to create values for the
parameters $k_n$ in the respective intervals in Eq.~(\ref{k1k2k3k4}).
Proceeding in this way, for $N=4$ ($p=2$), we obtain the values $k_1=1.00725$,
$k_2=1.04665$, $k_3=1.12192$, and $k_4=1.3129$, and a precision of $0.00001$,
which is, again, the lowest precision for a given $N$. As becomes obvious,
the lowest precisions are obtained for uniform partitions. This is not only
the case for $N$ being a power of $2$ but also for $N$ being a power of $3$, as
discussed in the next section.

\begin{figure}[ht]
\includegraphics[width=0.7\linewidth]{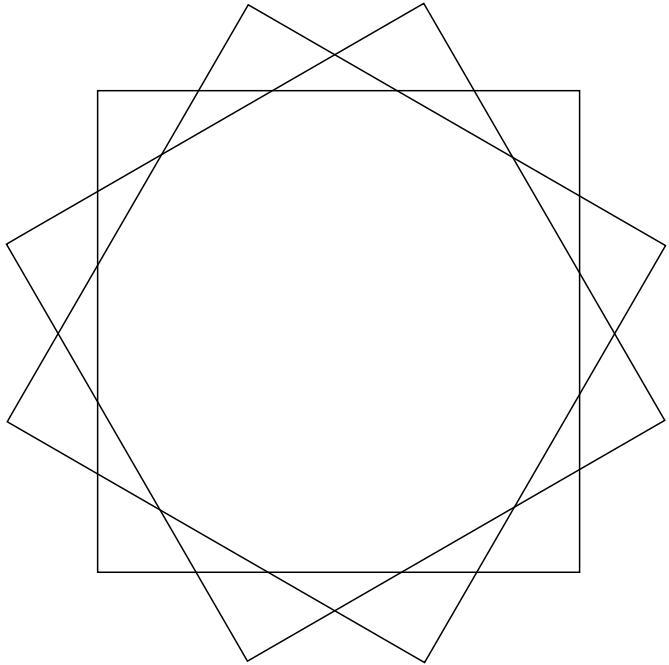}
\caption{The case of three ``boxes'' and accordingly three approximating
  exponents. The picture becomes more and more rotationally symmetric as the
  number of boxes grows.\label{fig2}}
\end{figure}

\section{A similar ternary procedure}
As the approximation~(\ref{rteq3}) gained high precision, we tried and
succeeded in finding a geometric interpretation for this, as is shown in
Fig.~\ref{fig2}. In this ternary approach, the initial step is to rotate the
square not by an angle of $\pi/4$ as in the previous approach but by an
angle of $\pi/6$. The overlapping squares in Fig.~\ref{fig2} can be split up
into three 12-angle figures that, at the same time, determine the constraints
for the parameters $k_i$,
\begin{equation}
1\le k_1\le 1/\cos(\pi/12)\le k_2\le 1/\cos(\pi/6)\le k_3\le 1/\cos(\pi/4).
\end{equation}
Note that the values $k_1=1.02335$, $k_2=1.05674$ and $k_3=1.28633$ chosen in
Eq.~(\ref{rteq3}) fit into these intervals. This procedure can be continued
iteratively in a ternary way,  i.e., turning the 12-angle figures by an angle
of $\pi/18$, and generally by the angle $\alpha=\pi/(2\cdot 3^p)$. In the
next section, we deal with the convergence of this and the previous procedure
for increasing values of $p$.

\section{Conclusions and Outlook}
In this paper, we have given an approximation for the Gauss integral with a
finite boundary in terms of the square root of a normalized sum of normal
distributions plus one, each of those distributions depending on the
(symmetric) boundary $[-t,t]$ of the integral and a set of maximally $N$
parameters $k_n$. By simple geometrical means, it is shown that these
parameters are constrained to intervals given by the inverse cosine with
equally distributed angles. We performed this approximation procedure in
both a binary ($N=2^p$) and a ternary way ($N=3^p$) and showed that the
procedure converges for an increasing degree $p$. The continuum limit leads
to a further approximation.

\subsection*{Acknowledgments}
The research was supported in part by the European Regional Development Fund
under Grant No.~TK133.

\begin{appendix}
\section{On the convergence of the procedure}
In general, one has
\begin{equation}
P^2(k_1,k_2,\ldots,k_N,t)=1-\frac{1}{N}\,\sum^{N}_{n=1}e^{-k_n^2t^2/2}\,.
\label{Eqf2}
\end{equation}
For geometry reasons, in the limit $N\to \infty$, the largest parameter in the
infinite set $\{k_1,k_2,\ldots,k_{\infty}\}$ must be $\sqrt{2}$, whereas the
lowest one must be $1$. The reason is that, in using the technique as in
Fig.~\ref{fig1}(B) over and over again, the final areas of integration turn to
perfect circles between the radii $t$ and $t\sqrt{2}$. The convergence of this
method becomes obvious by considering the backstep iteration. Suppose we start
with an approximation for a given set of parameters $k_n$ with a given
precision. The degeneration of two adjacent parameters means that a partition
is skipped, leading to a more imprecise approximation as the degree of freedom
in choosing different parameter values is lost.

For the general analysis we calculate the convergence by fixing the parameters
in Eq.~(\ref{Eqf2}) to the upper boundary, $k_n=1/\cos(\pi n/(4\cdot 2^p))$,
and analyse
\begin{equation}
\Delta_N(N,t)=\Delta(k_1,\ldots,k_N,t)=P(k_1,\ldots,k_N,t)-P(t)
\end{equation}
for $N=2^p$ and a fixed value of $t$, e.g., $t_0=1.0668$, at which, the
uncertainty range of Eq.~(\ref{Eq44}) turns out to be maximal. One obtains the
values in Table~\ref{tab1} demonstrating the convergence of the approximations.

\begin{table}[ht]
\caption{\label{tab1}Deviations in the uniform $N=2^p$ approximations for
increasing $p$. Note that values higher than $p=15$ could not be checked
with the PC at hand.}
\begin{tabular}{c|ccccc}\hline
\textbf{$p$}&$11$&$12$&$13$&$14$&$15$\\\hline
\textbf{$|\Delta_N(2^p,t_0)|$}&$0.00004$&$0.00002$&$0.00001$
  &$0.000005$&$0.0000026$\\\hline
\end{tabular}
\end{table}

For the ternary procedure we again fix the parameters to the upper
boundary, $k_n=1/\cos(\pi n/(4\cdot 3^p))$, and analyse $\Delta_N(N,t)$ for
$N=3^p$ and for the same fixed value $t_0=1.0668$. One obtains the values in
Table~\ref{tab2}.

\begin{table}[ht]
\caption{Deviations in the uniform $N=3^p$ approximations for
increasing $p$. Note that values higher than $p=10$ could not be checked
with the PC at hand.\label{tab2}}
\begin{tabular}{c|ccccccc}\hline
\textbf{$p$}&$6$&$7$&$8$&$9$&$10$\\\hline
\textbf{$|\Delta_N(3^p,t_0)|$}&$0.0001$&$0.00004$&$0.00001$
  &$0.000004$&$0.000001$\\\hline
\end{tabular}
\end{table}

The values in Tables~\ref{tab1} and~\ref{tab2} can be approximated by the
formula $|\Delta_p(N,1)|<0.09/N$, i.e., the deviation is inversely proportional
to $N$. This can be seen as follows. The worst error of the squared Gauss
integral $P^2(t)$ is given by using the approximations where the $k_n$ takes
the maximal or minimal values, respectively. The difference between these
squares of extremal values is given by
\begin{equation}
P_M^2(t)-P_m^2(t)=\frac1{N}\sum_{n=1}^N\left(e^{-(k_n^{\rm min})^2t^2/2}
  -e^{-(k_n^{\rm max})^2t^2/2}\right)=\frac{H(t)}N,
\end{equation}
where $H(t)=e^{-t^2/2}-e^{-t^2}$ is obtained by using the property
$k_n^{\rm min}=k_{n-1}^{\rm max}$ in order to cancel intermediate consecutive
terms. Using the third binomial to obtain
$P_M^2(t)-P_m^2(t)\approx 2P(t)(P_M(t)-P_m(t))$, one has
\begin{equation}
P_M(t)-P_m(t)\approx\frac{H(t)}{2P(t)N}.
\end{equation}
Finally, one can use $P(t)>P_{\rm min}(t)=\sqrt{1-e^{-t^2/2}}$ to obtain
\begin{equation}
P_M(t)-P_m(t)<\frac1{2N}H(t)\sqrt{1-e^{-t^2/2}}<\frac{0.09}N,
\end{equation}
where $H(t)\sqrt{1-e^{-t^2/2}}\le\sqrt{2^23^3/5^5}$ is used. The error for
$P(t)$ itself is, at most, the difference between the two extremal values.
Eq.~(\ref{Eqf2}) can be considered as the discretised form of the Gauss
integral. Applying the continuum limit
$\sum_if_i(z_i)\Delta z_i\to\int f(z)dz$, one obtains
\begin{equation}
P^2(t)=1-\frac4\pi\int_0^{\pi/4}\exp\pfrac{-t^2}{2\cos^2\phi}d\phi.
\end{equation}
The exponential function can be expanded into a series of finite degree $N$.
Again, we obtain an approximation, as
\begin{equation}
\Big|P^2(t)-\frac4\pi\sum_{n=1}^N\frac{(-1)^{n-1}}{n!}\pfrac{t^2}2^nc_n\Big|
  <\frac{t^{2N}}{N!N}
\end{equation}
with
\begin{equation}
c_n=\int_0^{\pi/4}\frac{d\phi}{\cos^{2n}\phi}\ =\ \sum_{k=0}^{n-1}\frac1{2k+1}
  \begin{pmatrix}n-1\\ k\\\end{pmatrix}\ =\ {}_2F_1(1/2,1-n;3/2;-1),
\end{equation}
where ${}_2F_1(a,b;c;z)$ is the hypergeometric function.

\end{appendix}

\end{document}